\documentclass[12pt]{article}

\usepackage{amssymb}
\font\smallit=cmti10
\font\smalltt=cmtt10

\def\B{\hfill $\Box$}

\begin{document}
\begin{center}
{\uppercase{\bf Two Color Off-diagonal Rado-type Numbers}}

\vskip 20pt
{\bf Kellen Myers\footnote{
This work was done as part of a summer REU,
funded by Colgate University,
while the first author was an undergraduate at Colgate University,
under the directorship of the second author.
}}\\
\vskip 5pt
and
\vskip -5pt
{\bf Aaron Robertson}\\
{\smallit Department of Mathematics,
Colgate University,
Hamilton, NY 13346}\\
{\smalltt aaron@math.colgate.edu}
\end{center}

\begin{abstract}
We show that for any two linear homogenous equations $\mathcal{E}_0,\mathcal{E}_1$, each with at least three
variables and coefficients not all the same sign, any 2-coloring of $\mathbb{Z}^+$ admits monochromatic
solutions of color 0 to $\mathcal{E}_0$ or monochromatic solutions of color 1 to $\mathcal{E}_1$. We define
the 2-color off-diagonal Rado number $RR(\mathcal{E}_0,\mathcal{E}_1)$ to be the smallest $N$ such that
$[1,N]$ must admit such solutions. We determine a lower bound for $RR(\mathcal{E}_0,\mathcal{E}_1)$ 
in certain cases when
each $\mathcal{E}_i$ is of the form $a_1x_1+\dots+a_nx_n=z$ as well as find the exact value of
$RR(\mathcal{E}_0,\mathcal{E}_1)$ when each is of the form $x_1+a_2x_2+\dots+a_nx_n=z$. We then present a
Maple package that determines upper bounds for off-diagonal Rado numbers of a few particular types, and use
it to quickly prove two previous results for diagonal Rado numbers.

\end{abstract}
\vskip 30pt
\section*{\normalsize 0. Introduction}

For $r \geq 2$, an $r$-coloring of the positive integers $\mathbb{Z}^+$ is an assignment 
$\chi : \mathbb{Z}^+ \rightarrow \{0,1,\dots,r-1\}$. Given a diophantine equation $\mathcal{E}$
in the variables $x_1,\dots,x_n$, we say a solution
$\{\bar{x}_i\}_{i=1}^n$ is monochromatic if $\chi(\bar{x}_i)=\chi(\bar{x}_j)$ for every $i,j$ pair. A
well-known theorem of Rado states that, for any $r \geq 2$, a linear homogeneous equation
$c_1x_1+\dots+c_nx_n=0$ with each
$c_i\in
\mathbb{Z}$ admits a monochromatic solution in $\mathbb{Z}^+$
 under any $r$-coloring of $\mathbb{Z}^+$ if and only if some nonempty subset of $\{c_i\}_{i=1}^n$ sums
to zero. The smallest
$N$ such that any $r$-coloring of $\{1,2,\dots,N\}=[1,N]$ satisfies this condition is called the 
$r$-color Rado number
for the equation
$\mathcal{E}$. However, Rado also proved the following,
much lesser known, result.

\noindent{\bf Theorem 0.1} ({\it Rado [6]}) Let $\mathcal{E}$
 be a linear homogeneous equation
with integer coefficients. Assume that
$\mathcal{E}$ has at least 3 variables
with both positive and negative coefficients.
Then any $2$-coloring of $\mathbb{Z}^+$ admits
a monochromatic solution to $\mathcal{E}$.

\noindent
{\it Remark.}  Theorem 0.1 cannot be extended to more than $2$ colors, without restriction
on the equation.  For example,  Fox and Radoi\v ci\'{c} [2] have shown, in particular, that 
there exists a 3-coloring of $\mathbb{Z}^+$ that admits no monochromatic solution to
$x+2y=4z$.  For more information about equations that have
finite colorings of $\mathbb{Z}^+$ with no monochromatic solution see [1] and [2].

In [4], the 2-color Rado numbers are determined for equations of the form $a_1x_1+\dots+a_nx_n=z$
where one of the $a_i$'s is $1$.  The case when $\min(a_1,\dots,a_n)= 2$ is done in [5],
while the general case is settled in [3].

In this article, 
we investigate the ``off-diagonal" situation.
To this end, for $r \in \mathbb{Z}^+$ define an off-diagonal Rado number for the equations
$\mathcal{E}_i$,
$0\leq i\leq r-1$, to be the least integer
$N$ (if it exists) for which any $r$-coloring of $[1,N]$ must admit a monochromatic solution to 
$\mathcal{E}_i$ of color $i$ for some $i \in [0,r-1]$. In this paper, when $r=2$ we will prove the existence
of such numbers  and determine
particular values and lower bounds in several specific cases when the two equations are of the form
$a_1x_1+\dots+a_nx_n=z$.

\section*{\normalsize 1. Existence}

The authors were unable to find an English translation
of Theorem 0.1.  For the sake of completeness,
we offer a simplified version of Rado's original
proof.

\noindent
{\it Proof of Theorem 0.1} (due to Rado [6])
Let $\sum_{i=1}^k \alpha_i x_i = \sum_{i=1}^{\ell} \beta_i y_i$
be our equation, where $k \geq 2$, $\ell \geq 1$,
$\alpha_i \in \mathbb{Z}^+$ for $1 \leq i \leq k$,
and $\beta_i \in \mathbb{Z}^+$ for $1 \leq i \leq \ell$.
By setting $x=x_1=x_2=\cdots=x_{k-1}$, $y=x_{k}$, and
$z=y_1=y_2=\cdots=y_{\ell}$, we  may consider solutions to
$$
ax+by=cz,
$$
where $a=\sum_{i=1}^{k-1} \alpha_i$, $b=c_k$, and
$c = \sum_{i=1}^{\ell} \beta_i$.  We will
denote $ax+by=cz$ by $\mathcal{E}$.

Let $m = \mathrm{lcm}\left(\frac{a}{\gcd(a,b)},
\frac{c}{\gcd(b,c)}\right)$.  Let $(x_0,y_0,z_0)$
be the solution to $\,\mathcal{E}$ 
with $\max(x,y,z)$ a minimum, where the maximum
is taken over all solutions of positive integers
to $\mathcal{E}$.  Let $A=\max(x_0,y_0,z_0)$.

Assume, for a contradiction, that there exists
a $2$-coloring of $\mathbb{Z}^+$ with no
monochromatic solution to $\mathcal{E}$.
First, note that for any $n \in \mathbb{Z}^+$, 
the set $\{in: i=1,2,\dots,A\}$ cannot be monochromatic,
for otherwise $x=x_0n$, $y=y_0n$, and $z=z_0n$ is a monochromatic
solution, a contradiction.

Let $x=m$ so that $\frac{bx}{a}, \frac{bx}{c} \in \mathbb{Z}^+$.
Letting red and blue be our two colors, we may assume,
without loss of generality, that $x$ is red.
Let $y$ be the smallest number in
$\{im: i=1,2,\dots,A\}$ that is blue.
Say $y =\ell m$ so that $2 \leq  \ell \leq A$.

For some $n \in \mathbb{Z}^+$, we have that
 $z = \frac{b}{a} (y-x)n$ is
blue, otherwise $\{i \frac{b}{a}(y-x): i=1,2,\dots\}$
would be red, admitting a monochromatic
solution to $\mathcal{E}$. Then $w =
\frac{a}{c}z +
\frac{b}{c} y$ must be red, for otherwise $az+by=cw$ and
$z,y,$ and $w$ are all blue, a contradiction.
Since $x$ and $w$ are both red, we have that
$
q=\frac{c}{a}w - \frac{b}{a}x = \frac{b}{a} (y-x)(n+1)
$
must be blue, for otherwise $x,w,$ and $q$ give
a red solution to $\mathcal{E}$.
As a consequence, we see that
$
\left\{i\frac{b}{a}(y-x) : i=n,n+1,\dots\right\}
$
is monochromatic.  This gives us that
$
\left\{i\frac{b}{a} (y-x) n: i=1,2,\dots,A\right\}
$
is monochromatic, a contradiction.
\hfill{$\Box$}

Using the above result, we offer an ``off-diagonal"
consequence.

\noindent
{\bf Theorem 1.1} Let $\mathcal{E}_0$ and
$\mathcal{E}_1$ be linear homogeneous equations
with integer coefficients. Assume that
$\mathcal{E}_0$ and
$\mathcal{E}_1$ each have at least 3 variables
with both positive and negative coefficients.
Then any $2$-coloring of $\mathbb{Z}^+$ admits
either a solution to $\mathcal{E}_0$ of
the first color or a solution to $\mathcal{E}_1$
of the second color.

\noindent
{\it Proof.} 
Let $a_0,a_1,b_0,b_1,c \in \mathbb{Z}^+$ and
denote by $\mathcal{G}_i$ the equation
$a_ix+b_iy=cz$ for $i=0,1$.
 Via the same argument given in
the proof to Theorem 0.1, we may consider solutions
to $\mathcal{G}_0$ and $\mathcal{G}_1$.  (The coefficients
on $z$ may be taken to be the same in both equations by finding
the lcm of the original coefficients on $z$ and
adjusting the other coefficients accordingly.)

Let the colors be red and blue.  We want to show that
any $2$-coloring admits either a red solution to $\mathcal{G}_0$
or a blue solution to $\mathcal{G}_1$.  From Theorem 0.1,
we have monochromatic solutions to each of these equations.
Hence, we assume, for a contradiction, that any monochromatic
solution to $\mathcal{G}_0$ is blue and that any monochromatic
solution to $\mathcal{G}_1$ is red.  This gives us that
for any $i \in \mathbb{Z}^+$,
if $ci$ is blue, then $(a_1+b_1)i$ is red
(else we have a blue solution to $\mathcal{G}_1$).

Now consider monochromatic solutions in $c \mathbb{Z}^+$.
Via the obvious bijection between colorings of $c \mathbb{Z}^+$
and $\mathbb{Z}^+$ and the fact that
linear homogeneous equations are unaffected by
dilation, Theorem 0.1 gives us the existence of
monochromatic solutions in $c \mathbb{Z}^+$.
If $cx,cy,cz$ solve $\mathcal{G}_0$ and are the same
color, then they must be blue.  Hence,
$\hat{x}=(a_1+b_1)x, \hat{y}=(a_1+b_1)y,$ and 
$\hat{z}=(a_1+b_1)z$ are
all red.  But, $\hat{x},\hat{y},\hat{z}$ solve
$\mathcal{G}_0$.  Thus, we have a red solution
to $\mathcal{G}_0$, a contradiction.
\hfill{$\Box$}

\section*{\normalsize 2. Two Lower Bounds}

Given the results in the previous section, we make
a definition, which uses the following
notation.

\noindent
{\bf Notation} For $n \in \mathbb{Z}^+$ and
$\vec{a} = (a_1,a_2,\dots,a_n) \in \mathbb{Z}^n$,
denote by
$\mathcal{E}_n(\vec{a})$ the linear homogeneous equation
$
\sum_{i=1}^n a_i x_i = 0.
$

\noindent
{\bf Definition} 
For $k,\ell \geq 3, \vec{b} \in \mathbb{Z}^k,$
and $ \vec{c} \in \mathbb{Z}^\ell,$
we let $RR(\mathcal{E}_k(\vec{b}),\mathcal{E}_\ell(\vec{c}))$
be
the minimum integer $N$, if it exists, such that any $2$-coloring of
$[1,N]$ admits either a solution to $\mathcal{E}_k(\vec{b})$ of
the first color or a solution to $\mathcal{E}_\ell(\vec{c})$
of the second color.

We now develop a general lower bound for certain
types of
those numbers guarenteed to exist by Theorem 1.1.

\noindent
{\bf Theorem 2.1} For $k,\ell \geq 2$, let 
$ b_1,b_2,\dots,b_{k-1},c_1,c_2,\dots,c_{\ell-1} \in \mathbb{Z}^+.$ 
Consider
$\mathcal{E}_k=\mathcal{E}_k(b_1,b_2,$ $\dots,$ $b_{k-1},-1)$ and 
$\mathcal{E}_\ell=\mathcal{E}_\ell(c_1,c_2,\dots,c_{\ell-1},-1)$,
written 
so that $b_1 = \min(b_1,b_2,\dots,b_{k-1})$
and $c_1 = \min(c_1,c_2,\dots,c_{\ell-1})$.
Assume that $t=b_1=c_1$.
Let $q = \sum_{i=2}^{k-1} b_i$ and
$s =\sum_{i=2}^{\ell-1}c_i$.  
Let (without loss of generality) $q \geq s$.  
Then
$$
RR(\mathcal{E}_k,\mathcal{E}_\ell) \geq
t(t+q)(t+s)+s.
$$

\noindent
{\it Proof.}  
Let $N=t(t+q)(t+s)+s$ and
consider the  $2$-coloring of $[1,N-1]$
defined by coloring $[s+t, (q+t)(s+t)-1]$
red and its complement blue.  We will show that this
coloring avoids red solutions to $\mathcal{E}_k$ and
blue solutions to $\mathcal{E}_\ell$.

We first consider any possible red solution to $\mathcal{E}_k$.
 The value of $x_k$ would have to be at least $t(s+t)+q(s+t)=(q+t)(s+t)$.
Thus, there is no suitable red solution.
Next, we consider $\mathcal{E}_\ell$.  If $\{x_1,x_2,\dots,x_{\ell-1}\}
\subseteq [1,s+t-1]$, then
$x_\ell < (q+t)(s+t)$.  Hence, the
smallest possible blue solution to $\mathcal{E}_\ell$ 
has
$x_i \in [(q+t)(s+t),N-1]$ for some $i \in [1,\ell-1]$.  However, this
gives $x_\ell \geq t(q+t)(s+t)+s > N-1$.  Thus, there is no suitable
blue solution. 
\hfill $\Box$

The case when $k = \ell = 2$ in Theorem 2.1 can
be improved somewhat in certain cases,
depending upon the relationship between $t$, $q$,
and $s$.  This result is presented below.

\noindent
{\bf Theorem 2.2} 
Let $t,j \in \mathbb{Z}^+$.
Let
$\mathcal{F}^t_j$ represent the equation $tx+jy=z$. 
Let $q,s \in \mathbb{Z}^+$ with $q \geq s \geq t$. 
Define $m=\frac{\gcd(t,q)}{\gcd(t,q,s)}$. 
Then
$$
RR(\mathcal{F}^t_q,\mathcal{F}^t_s) \geq
t(t+q)(t+s)+ms.
$$

\noindent
{\it Proof.}  
Let $N=
t(t+q)(t+s)+ms$ and
consider the  $2$-coloring $\chi$ of $[1,N-1]$ 
defined by
coloring
$$
R=[s+t, (q+t)(s+t)-1] \cup \{t(t+q)(t+s)+is: 1 \leq i \leq m-1\}
$$
red and
$B=[1,N-1] \setminus R$ blue.
We will show that this
coloring avoids red solutions to $\mathcal{F}^t_q$ and
blue solutions to $\mathcal{F}^t_s$.

We first consider any possible red solution to $\mathcal{F}^t_q$.
 The value of $z$ would have to be at least $t(s+t)+q(s+t)=(q+t)(s+t)$
and congruent to $0$ modulo $m$.
Since $t(t+q)(t+s) \equiv 0 \, (\mbox{mod }m)$ but $is \not \equiv  0 \,(\mbox{mod }m)$ for
$1 \leq i \leq m-1$,
there is no suitable red solution. 
Next, we consider $\mathcal{F}^t_s$.  If $\{x,y\}
\subseteq [1,s+t-1]$, then
$s+t \leq z < (q+t)(s+t)$.  Hence, the
smallest possible blue solution to $\mathcal{F}^t_s$ has
$x$ or $y$ in $[(q+t)(s+t),N-1]$.  However, this
gives $z \geq t(q+t)(s+t)+s > N-1$.  
By the definition of the coloring,
$z$ must be red.
Thus, there is no suitable blue
solution to $\mathcal{F}^t_s$. 
\hfill $\Box$

\section*{\normalsize 3.  Some Exact Numbers}

In this section, we will determine some of the values
of $RR_1(q,s)=RR(x+qy=z,x+sy=z)$, where $1 \leq s \leq q$.
The subscript $1$ is present to emphasize
the fact that we are using $t=1$ as defined in
Theorem 2.1.  In this section we will
let $RR_t(q,s)=RR(tx+qy=z,tx+sy=z)$ and we will
denote the equation $tx+jy=z$ by $\mathcal{F}^t_j$.

\noindent
{\bf Theorem 3.1}  Let $1 \leq s \leq q$.  Then
$$ RR_1(q,s) = \left\{
\renewcommand{\arraystretch}{1.5}
\begin{array}{ll}
 2q +2 \left\lfloor \frac{q+1}{2}\right\rfloor + 1&\mbox{ for }s=1\\
 (q+1)(s+1)+s&\mbox{ for } s \geq 2. \\
\end{array}
\right.
$$

\noindent
{\it Proof.}  We start with the case $s=1$.
Let $N= 2q +2 \left\lfloor
\frac{q+1}{2}\right\rfloor + 1$.  
We first improve the lower bound given by Theorem 2.1
for this case.  

Let $\gamma$ be the 2-coloring of $[1,N-1]$ defined
as follows.  The first $2 \lfloor \frac{q+1}{2}\rfloor -1$
integers alternate colors with the color of $1$ being blue.
We then color $\left[2 \lfloor \frac{q+1}{2}\rfloor,2q+1\right]$
red.  We color the last $2 \lfloor \frac{q+1}{2}\rfloor -1$
integers with alternating colors, where the color of
$2q+2$ is blue. 

 First consider possible blue solutions to
$x+y=z$.  If $x,y \leq 2 \lfloor \frac{q+1}{2}\rfloor -1$,
then $z \leq 2q$.
Under $\gamma$, such a $z$ must be red.  Now, if
exactly one of $x$ and $y$ is greater than $2q+1$, then
$z$ is odd and greater than $2q+1$.  Again, such a $z$
must be red.  Finally, if both $x$ and $y$ are greater
than $2q+1$, then $z$ is too big.  Hence, $\gamma$ admits
no blue solution to $x+y=z$.

Next, we consider possible red solutions to $x+qy=z$.
If $x,y \leq  \lfloor \frac{q+1}{2}\rfloor -1$, then
$z$ must be even.  Also, since $x$ and $y$ must both
be at least $2$ under $\gamma$, we see that
$z \geq 2q+2$.  Under $\gamma$, such a $z$ must be
blue.  If one (or both) of $x$ or $y$ is greater than
$\lfloor \frac{q+1}{2}\rfloor -1$, then
$z \geq N-1$, with equality possible.  However, with
equality, the color of $z$ is blue.  Hence, $\gamma$
admits no red solution to $x+qy=z$.

We move onto the upper bound.
Let $\chi$ be a $2$-coloring
of $[1,N]$ using the colors red and blue.  Assume, for
a contradiction, that there is no red solution to
$\mathcal{F}^1_q$ and no blue solution to $\mathcal{F}^1_1$.    
We break the
argument into $3$ cases.

\noindent
{\tt Case 1.}  $1$ is red.  Then $q+1$ must be blue since
otherwise $(x,y,z)=(1,1,q+1)$ would be a red solution to
$\mathcal{F}^1_q$.  Since $(q+1,q+1,2q+2)$ satisfies $\mathcal{F}_1^1$,
we have that $2q+2$ must be red.  Now, since
$(q+2,1,2q+2)$ satisfies $\mathcal{F}^1_q$, we see that $q+2$ must be blue.
Since $(2,q+2,q+4)$ satisfies $\mathcal{F}_1^1$ we have that
$q+4$ must be red.  This implies that $4$ must be blue
since $(4,1,q+4)$ satisfies $\mathcal{F}_q^1$.
But then $(2,2,4)$ is a blue solution to $\mathcal{F}_1^1$,
a contradiction.

\noindent
{\tt Case 2.} $1$ is blue and $q$ is odd. 
Note that in this case we have $N=3q+2$.
Since $1$ is blue, $2$ must be
red, which, in turn, implies that $2q+2$ must be blue.
Since $(q+1,q+1,2q+2)$ solves $\mathcal{F}_1^1$,
we see that $q+1$ must be red.
Now, since $(j,2q+2,2q+j+2)$ solves $\mathcal{F}_1^1$
and $(j+2,2,2q+j+2)$ solves $\mathcal{F}_q^1$, we have
that for any $j \in \{1,3,5,\dots,q\}$,
the color of $j$ is blue.  With $2$ and $q$ both
red, we have that $3q$ is blue, which implies that
$3q+1$ must be red.  Since $(q+1,2,3q+1)$ solves $\mathcal{F}_q^1$,
we see that $q+1$ must be blue, and hence $q+2$ is red.
Considering $(q+2,2,3q+2)$, which solves $\mathcal{F}_q^1$,
and $(q,2q+2,3q+2)$, which solves $\mathcal{F}_1^1$,
we have an undesired monochromatic solution, a contradiction.

\noindent
{\tt Case 3.} $1$ is blue and $q$ is even. 
Note that in this case we have $N=3q+1$.  As in
Case 2, we argue that for any $j \in \{1,3,5,\dots,q-1\}$,
the color of $j$ is blue.  As in Case 2,
both $2$ and $q+1$ must be red, so that $3q+1$ must be blue.
But $(q-1,2q+2,3q+1)$ is then a blue solution to $\mathcal{F}_1^1$,
a contradiction.

Next, consider the cases when $s \geq 2$.  
From Theorem 2.1, we have
$RR_1(q,s) \geq  (q+1)(s+1)+s$.
  We proceed by showing that $RR_1(q,s) \leq  (q+1)(s+1)+s$.

In the  case when $s=1$ we used an obvious
``forcing" argument.  As such, we have automated
the process in the Maple package {\tt SCHAAL},
available for download from the second author's
webpage\footnote{{\tt http://math.colgate.edu/$\sim$aaron}}.
The package is detailed in the next subsection, but
first we finish the proof.
Using {\tt SCHAAL} we
find the following (where we use the fact that
$s \geq 2$):

\noindent
1)  If $1$ is red, then the elements in
$\{s,q+s+1,qs+q+s+1\}$ must be both red and blue,
a contradiction.

\noindent
2)  If $1$ is blue and $s-1$ is red, then the elements
in $\{1,2,2q-1,2s+1,2q+1,2q+2s-1,2q+2s+1\}$ must be both red and blue,
a contradiction.

\noindent
3)  If $1$ and $s-1$ are both blue, the analysis is a bit
more involved.  First, by assuming $s \geq 2$ we find
that $2$ must be red and $s$ must be blue.  Hence,
we cannot have $s = 2$ or $s=3$, since if $s=2$ then
$2$ is both red and blue, and if $s=3$ then
since $s-1$ is blue, we again have that $2$ is both
red and blue.  Thus, we may assume that $s \geq 4$.
Using {\tt SCHAAL} with $s \geq 4$ now produces the result that the
elements in
$
\{4, s + 1,q + 1,2 s - 1, 2 s, q + 2 s + 1,3 s + 1, 5 q + 1,4 q + s+ 1,4 q + 2 s - 1,
4 q + 2 s, 4 q + 3 s + 1, 5 q + 2 s + 1,  q s - 3 q + 1,
 q s - 3 q + 2 s +1,
q s - 3 q + s - 1, q s + q + 1,   
  q s + q + s - 1,  q s + q + 2 s + 1\}
$ must be both red and blue, a contradiction.

This completes the proof of the theorem.
\hfill $\Box$

Using the above theorem, we offer the following corollary.

\noindent
{\bf Corollary 3.2} For $k,\ell \in \mathbb{Z}^+$,
let $a_1,\dots,a_k,b_1,\dots,b_\ell \in \mathbb{Z}^+$.  
Assume $\sum_{i=1}^k a_i \geq \sum_{i=1}^\ell b_i$.  Then
$$
RR_1(x+\sum_{i=1}^k a_iy_i=z,x+\sum_{i=1}^\ell b_iy_i=z)=
\left\{
\begin{array}{ll}
\displaystyle
2 \sum_{i=1}^k a_i + 2 \left\lfloor \frac{\sum_{i=1}^k a_i + 1}{2}\right\rfloor +1&
\mbox{for } \displaystyle\sum_{i=1}^\ell b_i = 1\\
\displaystyle
\left( \sum_{i=1}^k a_i + 1 \right)\!\!\! \left( \sum_{i=1}^\ell b_i +1\right) + \sum_{i=1}^\ell b_i&
\mbox{for } \displaystyle\sum_{i=1}^\ell b_i \geq 2.\\
\end{array}
\right. 
$$

\noindent
{\it Proof.}  
We start by proving that the coloring
given in the proof of Theorem 3.1 which provides the lower bound
 for the case $s=1$
also provides (with a slight modification) a lower bound for the case when $\sum_{i=1}^\ell b_i = 1$.
In this situation, we must show that the coloring
where the first $2 \lfloor \frac{\sum_{i=1}^k a_i+1}{2}\rfloor -1$
integers alternate colors with the color of $1$ being blue.
We then color $\left[2 \lfloor \frac{\sum_{i=1}^k a_i+1}{2}\rfloor,2\sum_{i=1}^k a_i+1\right]$
red.  We color the last $2 \lfloor \frac{\sum_{i=1}^k a_i+1}{2}\rfloor -1$
integers with alternating colors, where the color of
$2\sum_{i=1}^k a_i+2$ is blue.  An obvious parity argument shows that
there is no blue solution to $x+y=z$ (this is the case when
$\sum_{i=1}^\ell b_i = 1$) exists, so it remains to show that no
red solution to $x+\sum_{i=1}^k a_iy_i=z$ exists under this coloring.
Now, if $x$ and all the $y_i$'s are less than  $2 \lfloor \frac{\sum_{i=1}^k a_i+1}{2}\rfloor$,
then $z$ would be even and have value at least $2\sum_{i=1}^k a_i+2$.
This is not possible, so at least one of $x,y_1,\dots,y_k$ must have
value at least $2 \lfloor \frac{\sum_{i=1}^k a_i+1}{2}\rfloor$.
If $x \geq 2 \lfloor \frac{\sum_{i=1}^k a_i+1}{2}\rfloor$, then
$z \geq 2 \sum_{i=1}^k a_i + 2 \left\lfloor \frac{\sum_{i=1}^k a_i + 1}{2}\right\rfloor$.
Hence, either $z$ is blue or too big.
So, assume, without loss of generality,
that $y_1 \geq 2 \lfloor\frac{\sum_{i=1}^k a_i+1}{2}\rfloor$.
If $a_1 = 1$, then $z=x+y_1+\sum_{i=2}^k a_i y_i \geq 2 + 2 \lfloor\frac{\sum_{i=1}^k a_i+1}{2}\rfloor
+ 2\sum_{i=2}^k a_i = 2 \lfloor\frac{\sum_{i=1}^k a_i+1}{2}\rfloor + 2\sum_{i=1}^k a_i$
and again either $z$ is blue or too big.
If $a_1 \geq 2$ (and we may assume that
$k \geq 2$ so that $\sum_{i=1}^k a_i + 1 \geq 4$), then $z = x + a_1y_1 + \sum_{i=2}^k a_i y_i
> a_1 \cdot 2 \lfloor \frac{\sum_{i=1}^k a_i+1}{2}\rfloor + 2\sum_{i=2}^k a_iy_i
\geq 2(a_1 + \lfloor \frac{\sum_{i=1}^k a_i+1}{2}\rfloor) + 2\sum_{i=2}^k a_iy_i
= 2\lfloor \frac{\sum_{i=1}^k a_i+1}{2}\rfloor) + 2\sum_{i=1}^k a_iy_i$
and $z$ is too big.

Next, by coupling the above lower bound
with Theorem 2.1 (using $t=1$), it remains to prove that the
righthand sides of the theorem's equations serve
as upper bounds for $N=RR_1(x+\sum_{i=1}^k a_iy_i=z,x+\sum_{i=1}^\ell b_iy_i=z)$.
Letting $q=\sum_{i=1}^k a_i$ and $s = \sum_{i=1}^\ell b_i$, any
solution to $x+qy=z$ (resp., $x+sy=z$) is a solution to
$x+\sum_{i=1}^k a_i y_i$ (resp., $x+\sum_{i=1}^\ell b_iy_i=z$) by
letting all $y_i$'s equal $y$.  Hence,
$N \leq RR_1(q,s)$ and we are done.
\hfill$\Box$

\noindent
{\it Remark.}  When $a_i=1$ for $1 \leq i\leq k$, $\ell =1$, and
$b_1=1$ the numbers in Corollary 3.2 are
called the off-diagonal generalized Schur numbers.  In this case,
the values of the numbers have been determined [7].

\subsection*{\normalsize 3.1 About the Maple Package {\tt SCHAAL}}

This package is used to try to automatically provide
an upper bound for the off-diagonal Rado-type numbers
$RR_t(q,s)$.  The package employs a set of rules to follow,
while the overall approach is an implementation of
the above ``forcing" argument.

Let $t \geq 2$ be given, keep $q \geq s$
as parameters, and define $N = 
tqs+t^2q+(t^2+1)s+t^3$.
We let $\mathcal{R}$ and $\mathcal{B}$ be the set 
of red, respectively blue, elements in $[1,N]$.
The package {\tt SCHAAL} uses  the following rules.

\noindent
For $x,y \in \mathcal{R}$,

R1)  if $q | (y-tx)$ and $y-tx>0$, then $\frac{y-tx}{q} \in \mathcal{B}$; 

R2)  if $t | (y-qx)$ and $y-qx>0$, then $\frac{y-qx}{t} \in \mathcal{B}$;

R3)  if $(q+t) | x$ then $\frac{x}{q+t} \in \mathcal{B}$.

\noindent
For $x,y \in \mathcal{B}$,

B1)  if $s | (y-tx)$ and $y-tx>0$, then $\frac{y-tx}{s} \in \mathcal{R}$; 

B2)  if $t | (y-sx)$ and $y-sx>0$, then $\frac{y-sx}{t} \in \mathcal{R}$;

B3)  if $(s+t) | x$ then $\frac{x}{s+t} \in \mathcal{R}$.

We must, of course, make sure that the
 elements whose colors are implied
by the above rules are in $[1,N]$. This
is done by making sure that the
coefficients of $qs$, $q$, and $s$,
as well as the constant term are nonnegative
and at most equal to the corresponding coefficients
in $tqs+t^2q+(t^2+1)s+t^3$ (hence the need for
$t$ to be an integer and not a parameter).  
See the Maple code for
more details.

The main program of {\tt SCHAAL} is {\tt dan}.  The program
{\tt dan} runs until 
$\mathcal{R} \cap \mathcal{B} \neq \emptyset$ or
until none of the above rules
produce a color for a new element.

\section*{\normalsize 3.2 Some Diagonal Results Using {\tt SCHAAL}}

Included in the package {\tt SCHAAL} is the program
{\tt diagdan}, which is a cleaned-up version of {\tt dan}
in the case when $q=s$.  Using {\tt diagdan} we are able
to reprove the main results found in [4] and [5].
However, our program is not designed to reproduce
the results in [3], which keeps $t$ as a parameter and 
confirms the conjecture of Hopkins and Schaal [4]
that $R_t(q,q)=tq^2+(2t^2+1)q+t^3$.

\noindent
{\bf Theorem 3.3} ({\it Jones and Schaal [5]})
$R_1(q,q)= q^2+3q+1$

\noindent
{\it Proof.} By running {\tt diagdan}$(\{1\},\{\},1,q)$ we
find immediately that the elements in $\{1,2,q,2q+1,q^2+2q+1\}$
must be both red and blue, a contradiction.
\B

\noindent
{\bf Theorem 3.4} ({\it Hopkins and Schaal [4]})
$R_2(q,q) = 2q^2+9q+8$

\noindent
{\it Proof.} By running {\tt diagdan}$(\{1\},\{q\},2,q)$ we
find immediately that the elements in $\{q+2,2q^2+5q,\frac{1}{2}
(q^2+3q)\}$
must be both red and blue.  We then run
{\tt diagdan}$(\{1,q\},\{\},2,q)$ and find that
the elements in
$\{2,q+2,2q,6q,q^2+6q\}$
must be both red and blue.  The program ran for about 10
seconds to obtain this proof.
\B

\section*{\normalsize 3.3 Some Values of $RR_t(q,s)$}

We end this section (and paper) with some values of
$RR_t(q,s)$ for small values
of $t,q$ and $s$.

$$
\begin{array}{l|l|l|r|c|l|l|l|r}
t&q&s&\mbox{Value}&\hskip 30pt&t&q&s&\mbox{Value}\\
\hline 
2&3&2&43&&3&5&4&172\\
2&4&2&50&&3&6&4&201\\
2&5&2&58&&3&7&4&214\\
2&6&2&66&&3&8&4&235\\
2&7&2&74&&3&9&4&264\\
2&8&2&82&&3&10&4&277\\
2&9&2&90&&3&6&5&231\\
2&10&2&98&&3&7&5&245\\
2&4&3&66&&3&8&5&269\\
2&5&3&73&&3&9&5&303\\
2&6&3&86&&3&10&5&317\\
2&7&3&93&&3&7&6&276\\
2&8&3&106&&3&8&6&303\\
2&9&3&112&&3&9&6&330\\
2&10&3&126&&3&10&6&357\\
2&5&4&88&&3&8&7&337\\
2&6&4&100&&3&9&7&381\\
\end{array}
$$
\centerline{\bf Table 1:  Small Values of $RR_t(q,s)$}

$$
\begin{array}{l|l|l|r|c|l|l|l|r}
t&q&s&\mbox{Value}&\hskip 30pt&t&q&s&\mbox{Value}\\
\hline
2&7&4&112&&3&10&7&397\\
2&8&4&124&&3&9&8&420\\
2&9&4&136&&3&10&8&437\\
2&10&4&148&&3&10&9&477\\
2&6&5&122&&4&5&4&292\\
2&7&5&131&&4&6&4&324\\
2&8&5&150&&4&7&4&356\\
2&9&5&159&&4&8&4&388\\
2&10&5&178&&4&9&4&432^*\\
2&7&6&150&&4&10&4&452\\
2&8&6&166&&4&6&5&370\\
2&9&6&182&&4&7&5&401\\
2&10&6&198&&4&8&5&452\\ 
2&8&7&194&&4&9&5&473\\
2&9&7&205&&4&10&5&514\\
2&10&7&230&&4&7&6&446\\
2&9&8&228&&4&8&6&492\\
2&10&8&248&&4&9&6&526\\
2&10&9&282&&4&10&6&566\\
3&4&3&129&&4&8&7&556\\
3&5&3&147&&4&9&7&579\\
3&6&3&165&&4&10&7&630\\
3&7&3&192^*&&4&9&8&632\\
3&8&3&201&&4&10&8&680\\
3&9&3&219&&4&10&9&746\\
3&10&3&237&&5&11&5&820^*
\end{array}
$$
\normalsize
\centerline{\bf Table 1 cont'd:  Small Values of $RR_t(q,s)$}

\vskip 20pt
These values were calculated by matching Theorem 2.2's lower bound
with the Maple package {\tt SCHAAL}'s upper bound.
We use {\tt SCHAAL}
by letting $1$ be red and then letting $1$ be blue.
In many cases this is sufficient, however in many of the
remaining cases, we must consider subcases depending upon
whether $2$ is red or blue.  If this is still not sufficient,
we consider subsubcases depending upon whether the value in
the table, the integer $3$, the integer $4$, or the
integer $5$, is red or blue.  This is sufficient for
all values in Table 1, expect for those marked with an $^*$.
This is because,
except for those three values marked with an $^*$, all values agree
with the lower bound given by Theorem 2.2.  
For these three exceptional values, we can increase the
lower bound given in Theorem 2.2.

\noindent
{\bf Theorem 3.3}  Let $t \geq 3$.  Then
$
R_t(2t+1,t) \geq 6t^3+2t^2+4t.
$

\noindent
{\it Proof.}
It is easy to check that the $2$-coloring of
$[1,6t^3+2t^2+4t-1]$ defined by coloring
$\{1,2,6t\} \cup \{6t+3,\dots,6t^2+2t-1\} \cup
\{6t^2+2t \leq i \leq 12t^2+4t: i \equiv 0 \, (\mbox{mod } t)\}$
red and its complement blue avoids red solutions to
$tx+(2t+1)y=z$ and blue solutions to $tx+ty=z$.
(We use $t > 2$ so that $6t$ is the minimal red element
that is congruent to $0$ modulo $t$.)
\hfill $\Box$

\noindent
{\it Remark.}   The lower bound in the above theorem is not tight.
For example, when $t=6$, the $2$-coloring of $[1,1392]$
given by coloring $\{1,2,3,37,39,40,41,43,46,47,48,49,50,52,56\}
\cup [58,228] \cup \{234 \leq i \leq 558:  i \equiv 0\,(\mbox{mod }6)\}
\cup \{570,576,594,606,612,648,684\}$ red and its complement blue
avoids red solutions to $6x+13y=z$ and blue solutions to
$6x+6y=z$.  Hence, $RR_t(2t+1,t) > 6t^3+2t^2+4t$ for $t=6$.

We are unable to explain why $(b,c)=(2t+1,t)$
produces these ``anomolous" values while
others, e.g., $(b,c) = (2t-1,t)$, appear not
to do so.

\section*{\normalsize References}
\footnotesize
\parindent=0pt

[1] B. Alexeev, J. Fox, and R. Graham,
On Minimal colorings Without Monochromatic Solutions to
a Linear Equation, preprint.

[2] J. Fox and R. Radoi\v ci\'{c}, The Axiom of Choice and
the Degree of Regularity of Equations of the Reals,
preprint.

[3] S. Guo and Z-W. Sun, Determination of the
Two-Color Rado Number for
$a_1x_1+\cdots+a_mx_m=x_0$,
preprint, {\tt arXiv:math.CO/0601409}.

[4] B. Hopkins and D. Schaal, On Rado Numbers for
$\sum_{i=1}^{m-1} a_ix_i = x_m$, {\it Adv. Applied
Math.} {\bf 35} (2005), 433-441.

[5] S. Jones and D. Schaal, Some 2-color Rado Numbers,
{\it Congr. Numer.} {\bf 152} (2001), 197-199.

[6] R. Rado, Studien zur Kombinatorik, {\it Mathematische Zeitschrift}
{\bf 36} (1933), 424-480.

[7] A. Robertson and D. Schaal, Off-Diagonal Generalized Schur Numbers,
{\it Adv. Applied Math.} {\bf 26}, 252-257.

\end{document}